\documentclass[11pt]{article}

\usepackage{graphicx}
\DeclareGraphicsExtensions{.ps,.eps}

\usepackage{isolatin1}
\usepackage{a4}
\usepackage{latexsym}
\usepackage{amsfonts}
\usepackage{amsmath}
\usepackage{amsthm}
\usepackage{amssymb}
\usepackage{bbm}
\usepackage{verbatim}
\usepackage{algorithm}
\usepackage{algorithmic}
\usepackage{stmaryrd}
\usepackage{quantizerDesign}
\usepackage{quantizerDesignic}

\font\tenmath=msbm10 scaled 1200
\font\sevenmath=msbm7 scaled 1200
\font\fivemath=msbm5 scaled 1200

\newtheorem{Remark}{Remark}[section]

\newfam\mathfam \textfont\mathfam=\tenmath
\scriptfont\mathfam=\sevenmath \scriptscriptfont\mathfam=\fivemath
\def\math{\fam\mathfam}
\def\R{{\math R}}
\def\N{{\math N}}
\def\E{{\math E}}

\def\P{{\math P}}

\def \^#1{\if#1i{\accent"5E\i}\else{\accent"5E#1}\fi}

\def \lesim{\stackrel{<}{\sim}}

\def \ni{\noindent}

\providecommand{\norm}[1]{\lVert#1\rVert}

\providecommand{\biggnorm}[1]{\biggl\lVert#1\biggr\rVert}

\providecommand{\scalar}[2]{< #1, #2 >}

\providecommand{\scalarOn}[3]{< #1, #2 >_{#3}}
\providecommand{\indicator}[2]{\mathbbm{1}_{#1}(#2)}

\newcommand{\bigtimes}{\prod}
\newcommand{\Rd}{\R^d}

\newcommand{\Prob}{\mathbb{P}}

\newcommand{\Normal}{\mathcal{N}}
\newcommand{\ProdNormal}{\bigotimes_{j=1}^\infty \Normal(0, \lambda_j)}
\newcommand{\ai}{\alpha^{(i)}}

\DeclareMathOperator{\card}{card}

\DeclareMathOperator{\myspan}{span}

\renewcommand{\thequantizerDesign}{\Roman{quantizerDesign}}
\newcommand{\QD}{\thequantizerDesign}

\newtheorem{Thm}{Theorem}

\newtheorem{Lem}{Lemma}
\newtheorem{Pro}{Proposition}

\author{\sc
{\sc Harald Luschgy}\thanks{Universit\"at Trier, FB IV-Mathematik, D-54286 Trier, Germany.
E-mail: {\tt luschgy@uni-trier.de}}
\quad ,
\quad {\sc Gilles Pag\`es} \thanks{Laboratoire de Probabilit\'es et Mod\`eles al\'eatoires, UMR~7599, Universit\'e Paris 6, case 188, 4,
pl. Jussieu, F-75252 Paris Cedex 5. E-mail:{\tt  gpa@ccr.jussieu.fr}}
\quad {and} 
\quad {\sc Benedikt Wilbertz}\thanks{Universit\"at Trier, FB IV-Mathematik,
D-54286 Trier, Germany. E-mail: {\tt wilbertz@uni-trier.de}}
}

\title{\bf  Asymptotically optimal quantization schemes for Gaussian
processes\thanks{This work was supported in part by the AMaMeF Exchange Grant
1323 of the ESF.}}

\begin{document}

\maketitle

\bibliographystyle{plain}

\begin{abstract}
We describe quantization designs which lead to asymptotically and order optimal functional quantizers. Regular variation of
the eigenvalues of the covariance operator plays a crucial role to achieve these rates. 
For the development of a constructive quantization scheme we 
rely on the knowledge of the eigenvectors of the covariance operator in order to
transform the problem into a finite dimensional quantization problem of normal
distributions. 
\end{abstract}


\bigskip
\ni {\em MSC: 60G15, 60E99.}

\section{Introduction}
\setcounter{equation}{0}
\setcounter{Assumption}{0}
\setcounter{Theorem}{0}
\setcounter{Proposition}{0}
\setcounter{Corollary}{0}
\setcounter{Lemma}{0}
\setcounter{Definition}{0}
\setcounter{Remark}{0}

Functional quantization of stochastic processes can be seen as a discretization of the path-space of a process and the
 approximation (coding) of a process by finitely many deterministic functions from its path-space. In a Hilbert space setting this reads as follows.

Let $(H, < \cdot , \cdot >) $ be a separable Hilbert space with norm
$\|  \cdot \| $
and let $X : (\Omega, {\cal A}, \P) \rightarrow H$ be a random vector taking its values in $H$ with distribution $\P_X$. For $n \in \N$, the
$L^2$-quantization problem for $X$ of level $n$ (or of nat-level $\log n$) consists in minimizing
\[
\left(\E \min_{a \in \alpha} \| X - a \|^2 \right)^{1/2} =
\|
\min_{a \in \alpha} \| X-a \|  \|_{L^2(\P)}
\]
over all subsets $\alpha \subset H$ with $\mbox{card} (\alpha) \leq n$. Such a set $\alpha$ is called $n$-codebook or $n$-quantizer.
The minimal $n$th quantization
error of $X$ is then defined by
\begin{equation}
e_n (X) := \inf \left\{ ( \E \min_{a \in \alpha} \| X-a \|^2)^{1/2} :
\alpha \subset H, \; \mbox{card} (\alpha)
\leq n \right\} .
\end{equation}
Under the integrability condition
\begin{equation}\label{(1.2)}
\E \,\|  X \|^2 < \infty
\end{equation}
the quantity $e_n (X)$ is finite.

For a given $n$-quantizer $\alpha$ one defines an associated closest neighbour projection
\[
\pi_\alpha := \sum\limits_{a \in \alpha} a \mbox{\bf 1}_{C_{a} (\alpha)}
\]
and the induced $\alpha$-quantization (Voronoi quantization) of $X$ by
\begin{equation}
\hat{X}^\alpha := \pi_\alpha (X) ,
\end{equation}
where $\{ C_a (\alpha) : a \in \alpha\}$ is a Voronoi partition induced by $\alpha$, that is a Borel partition of $H$ satisfying
\begin{equation}
C_a (\alpha) \subset V_a (\alpha) := \{ x \in H : \| x-a \| = \min_{b \in \alpha} \| x-b \|  \}
\end{equation}
for every $a \in \alpha$. Then one easily checks that, for any random vector $X^{'} : \Omega \rightarrow \alpha \subset H$,
\[
\E\,\| X - X^{'} \|^2 \geq \,\E\,\| X - \hat {X}^\alpha \|^2 = \E\,\min_{a \in \alpha} \| X-a \|^2
\]
so that finally
\begin{eqnarray}
e_n (X) & = & \inf \left\{ ( \E\,\| X - \hat{X} \|^2)^{1/2} : \hat{X} = f(X), f : H \rightarrow H \;
\mbox{Borel measurable,} \; \right.\\
             && \left.\qquad \mbox{card} (f(H)) \leq n \right\} \nonumber \\
        & = & \inf \left\{ ( \E\,\| X- \hat{X} \|^2)^{1/2} : \hat{X} : \Omega \rightarrow H \;
\mbox{random vector,} \;
             \mbox{card} (\hat{X}(\Omega)) \leq n \right\} . \nonumber
\end{eqnarray}
Observe that the Voronoi cells $V_a (\alpha), a \in \alpha$ are closed and convex (where convexity is a characteristic
feature of the underlying Hilbert structure).
Note further that there are infinitely many $\alpha$-quantizations of $X$ which all produce the same quantization
error and $\hat{X}^\alpha$
is $\P$-a.s. uniquely defined if $\P_X$ vanishes on hyperplanes.

A typical setting for functional quantization is $H = L^2([0,1], dt)$ but is obviously not restricted to the Hilbert space setting.
Functional quantization is the natural extension to stochastic processes of the so-called optimal vector
quantization of random vectors in $H = \R^d$ which has been extensively investigated since the late 1940's in Signal
processing and Information Theory (see \cite{GERSH}, \cite{GRAY}).
For the mathematical aspects of vector quantization in $\R^d$, one may consult \cite{Foundations}, for algorithmic aspects see \cite{pagesOQ}
and ''non-classical'' applications can be found in \cite{pagesIntegration}, \cite{page2}.
For a first  promising application of functional quantization to the pricing of financial derivatives through numerical integration on
path-spaces see
\cite{FQFinance}. 

We address the issue of high-resolution quantization which concerns the performance of $n$-quantizers and the
behaviour of $e_n(X)$ as $n \rightarrow \infty$.
The asymptotics of $e_n(X)$ for $\R^d$-valued random vectors has been completely elucidated for non-singular distributions $\P_X$ by the Zador
Theorem (see \cite{Foundations}) and for a class of self-similar (singular) distributions by \cite{GRALU2}. In infinite dimensions no such global
results hold, even for Gaussian processes.

It is convenient to use the symbols $\sim$ and $\lesim$, where $a_n \sim b_n$ means $a_n/b_n \rightarrow 1$ and $a_n \lesim b_n$
means $\limsup_{n \to \infty} a_n/b_n \leq 1$. A measurable function $\varphi : (s, \infty) \rightarrow (0,
\infty)\, (s \geq 0)$ is said to be regularly varying at infinity with index $b \in \R$ if, for every $c > 0$,
\[
\lim_{x \to \infty} \frac{\varphi(cx) }{\varphi(x)} = c^b .
\]

Now let $X$ be centered Gaussian.
Denote by $K_X \subset H$ the reproducing kernel Hilbert space (Cameron-Martin space) associated to the covariance operator
\begin{equation}
C_{_X} : H \rightarrow H, \; C_{_X} y := \E\,(< \!y, X\! >\!X)
\end{equation}
of $X$. Let $\lambda_1 \geq \lambda_2 \geq \ldots > 0$ be the ordered nonzero eigenvalues of $C_{_X}$ and let $\{u_j : j \geq 1 \}$ be
the corresponding orthonormal basis of supp$(\P_X)$ consisting of eigenvectors (Karhunen-Lo\`eve basis). If
$d := \dim K_X < \infty$, then $e_n (X) = e_n \left( \bigotimes\limits^d_{j=1} \mathcal{N}(0, \lambda_j)\right)$,
the minimal $n$th $L^2$-quantization error of
$\bigotimes\limits^d_{j=1} \mathcal{N}(0, \lambda_j)$ with respect to the $l_2$-norm on $\R^d$, and thus we can read off the asymptotic
 behaviour of $e_n(X)$ from the high-resolution formula
\begin{equation}\label{1.7}
e_n (\bigotimes\limits^d_{j=1} \mathcal{N}(0, \lambda_j)) \sim q(d) \sqrt{2 \pi} \left( \Pi^{d}_{j=1}
\lambda_j\right)^{1/2 d}
\left( \frac{d+2}{d}\right)^{(d+2)/4} n^{-1/d} \; \mbox{ as } \; n \rightarrow \infty
\end{equation}
where $q(d) \in (0, \infty)$ is a constant depending only on the dimension $d$ (see \cite{Foundations}). Except in dimension $d=1$ and $d=2$, the
true value of $q(d)$ is unknown. However, one knows (see \cite{Foundations}) that
\begin{equation}\label{(1.8)}
q(d) \sim \left( \frac{d}{2 \pi e}\right)^{1/2} \; \mbox{ as } \; d \rightarrow \infty.
\end{equation}

Assume $\dim K_X = \infty$. Under regular behaviour of the eigenvalues the sharp asymptotics of $e_n (X)$ can be derived
analogously to (\ref{1.7}). In view of~(\ref{(1.8)}) it is reasonable to expect that the limiting constants can be evaluated. The recent high-resolution formula is as follows.

\begin{Thm}
(\cite{LUPA2}) Let $X$ be a centered Gaussian. Assume $\lambda_j \sim \varphi(j)$ as $j \rightarrow \infty$, where
$\varphi : (s, \infty) \rightarrow (0, \infty)$ is a decreasing, regularly varying function at infinity of
index $-b < -1$ for some $s \geq 0$. Set, for every $x > s$,
\[
\psi (x) := \frac{1}{x \varphi(x)} .
\]
Then
\[
e_n(X) \sim \left(\left( \frac{b}{2}\right)^{b-1}\!\! \frac{b}{b-1}\right)^{1/2} \psi ( \log n )^{-1/2}
\;
\mbox{ as }
\;
\vspace{0.5cm} n \rightarrow \infty .
\]
\end{Thm}

A high-resolution formula in case $b = 1$ is also available (see \cite{LUPA2}). Note that the restriction $-b \leq -1$
on the index of $\varphi$ is natural since $\sum\limits^\infty_{j=1}
\lambda_j < \infty$. The minimal $L^r$-quantization errors of $X$, $0 < r < \infty$, are strongly equivalent to the $L^2$-errors $e_n(X)$
(see \cite{DERE}) and thus exhibit the same high-resolution behaviour.

The paper is organized as follows. In Section~2 we investigate a quantization design,
 which furnishes asymptotically optimal quantizers in the situation of
Theorem 1. Here the Karhunen-Lo\`eve expansion plays a crucial
role. 
In section 3 we state different quantization designs,
which are all at least order-optimal
and discuss their possible implementations 
regarding the example of the Brownian motion.
The main focus in that section lies on "good" designs for finite $n\in\N$.

\section{Asymptotically optimal functional quantizers}
\setcounter{equation}{0}
\setcounter{Assumption}{0}
\setcounter{Theorem}{0}
\setcounter{Proposition}{0}
\setcounter{Corollary}{0}
\setcounter{Lemma}{0}
\setcounter{Definition}{0}
\setcounter{Remark}{0}
Let $X$ be a $H$-valued random vector satisfying~(\ref{(1.2)}). For every $n \in \N$, $L^2$-optimal $n$-quantizers $\alpha \subset H$ exist,
that is
\[
(\E\,\min_{a \in \alpha} \| X -a \|^2)^{1/2} = e_n (X) .
\]
If card (supp$(\P_X)) \geq n$, optimal $n$-quantizers $\alpha$ satisfy card$(\alpha) = n$,
$\P (X \in C_a (\alpha)) > 0$ and the stationarity condition
\begin{equation}\label{eq:KL}
a = \E\,(X \mid  \{ X \in C_a (\alpha)\}), \,a \in \alpha
\end{equation}
or what is the same
\begin{equation}\label{(3.1)}
\hat{X}^\alpha = \E\,(X \mid  \hat{X}^\alpha)
\end{equation}
for every Voronoi partition $\{ C_a(\alpha) : a \in \alpha \}$ (see \cite{LUPA1}). In particular,
$\E\,\hat{X}^\alpha = \E\,X$.

Now let $X$ be centered Gaussian with $\dim K_X = \infty$. The Karhunen-Lo\`eve basis $\{u_j : j \geq 1 \}$
consisting of normalized eigenvectors of $C_{_X}$ is optimal for the quantization of Gaussian random vectors
(see \cite{LUPA1}). So we start with the Karhunen-Lo\`eve expansion
\[
X \stackrel{H}{=} \sum\limits^\infty_{j=1} \lambda^{1/2}_j \xi_j u_j ,
\]
where $\xi_j = < \!X, u_j\!> / \lambda^{1/2}_j, j \geq 1$ are i.i.d. $\mathcal{N}(0,1)$-distributed random variables. The design of an
asymptotically optimal quantization of $X$ is based on optimal quantizing blocks of coefficients of variable ($n$-dependent)
block length. Let $n \in \N$ and fix temporarily
$m, l, n_1, \ldots , n_m \in \N$ with $\Pi^m_{j=1} n_j \leq n$, where $m$ denotes the number of blocks, $l$
the block length and $n_j$ the size of the quantizer for the $j$th block
\[
\xi^{(j)} := (\xi_{(j-1)l+1} , \ldots , \xi_{jl}), \quad j \in \{ 1, \ldots , m\}.
\]
Let $\alpha_j \subset \R^l$ be an $L^2$-optimal $n_j$-quantizer for
$\xi^{(j)}$ and let $\widehat{\xi^{(j)}} =  \widehat{\xi^{(j)}}^{\alpha_j}$
be a $\alpha_j$-quantization of $\xi^{(j)}$.
(Quantization of blocks $\xi^{(j)}$ instead of $( \lambda^{1/2}_{(j-1)l+1}
\xi_{(j-1)l+1} , \ldots , \lambda^{1/2}_{jl} \xi_{jl} )$ is asymptotically good
enough. 
For finite $n$ the quantization scheme will be considerably improved in Section
\ref{sec3}. ) Then, define a quantized version of $X$ by
\begin{equation}\label{3.2}
\hat{X}^n := \sum\limits^m_{j=1} \sum\limits^l_{k=1} \lambda^{1/2}_{(j-1)l+k} ( \widehat{\xi^{(j)}} )_k
u_{(j-1)l+k} .
\end{equation}
It is clear that card$(\hat{X}^n(\Omega )) \leq n$. Using (\ref{(3.1)}) for $\xi^{(j)}$, one gets
$\E\,\hat{X}^n = 0$. If
\[
\widehat{\xi^{(j)}} = \sum\limits_{b \in \alpha_j} b \mbox{\bf 1}_{C_b(\alpha_j)} (\xi^{(j)} ),
\]
then
\[
\hat{X}^n = \sum\limits_{a \in \times^m_{j=1} \alpha_j} ( \sum\limits^m_{j=1} \sum\limits^l_{k=1}
\lambda^{1/2}_{(j-1) l+k} a^{(j)}_k u_{(j-1) l+k} ) \Pi^m_{j=1} \mbox{\bf 1}_{C_{a^{(j)}} (\alpha_j)} (\xi^{(j)} )
\]
where $a = (a^{(1)}, \ldots , a^{(m)} ) \in \times^{m}_{j=1} \alpha_j$. Observe that in general,
$\hat{X}^n$ is not a Voronoi quantization of $X$ since it is based on the (less complicated) Voronoi
partitions for $\xi^{(j)}, j \leq m$. $(\hat{X}^n$ is a Voronoi quantization if $l = 1$ or if
$\lambda_{(j-1) l+1} = \ldots = \lambda_{jl}$ for every $j$.)
Using again (\ref{(3.1)}) for $\xi^{(j)}$ and the independence structure, one checks that $\hat{X}^n$
satisfies a kind of stationarity equation:
\[
\E\,(X \mid  \hat{X}^n) = \hat{X}^n .
\]
\begin{Lem} Let $n\ge1$. For every $l\ge 1$ and every $m\ge 1$
\begin{equation}\label{Pythagore}
\E\,\| X - \hat{X}^n \|^2 \leq \sum\limits^m_{j=1} \lambda_{(j-1)l+1} e_{n_{j}} (N(0, I_l))^2
+ \sum\limits_{j \geq ml+1} \lambda_j .
\end{equation}
Furthermore,~(\ref{Pythagore}) stands as an equality if $l=1$ (or $\lambda_{(j-1) l+1} = \ldots = \lambda_{jl}$ for every $j,\,l\ge 1$).
\end{Lem}

\noindent{\bf Proof.}
The claim follows from the orthonormality of the basis $\{ u_j : j \geq 1 \}$. We have
\[
\begin{array}{lcl}
\E\,\| X - \hat{X}^n \|^2
& = & \sum\limits^m_{j=1} \sum\limits^l_{k=1} \lambda_{(j-1) l+k} \E\,\mid \xi^{(j)}_k - (\widehat{\xi^{(j)}})_k
        \mid^2 + \sum\limits_{j \geq ml +1} \lambda_j \\
& \leq & \sum\limits^m_{j=1} \lambda_{(j-1)l+1} \sum\limits^l_{k=1} \E\,\mid \xi^{(j)}_k - \widehat{\xi^{(j)}} )_k \mid^2 +
        \sum\limits_{j \geq ml +1} \lambda_j \\
& = & \sum\limits^m_{j=1} \lambda_{(j-1)l+1} e_{n_{j}} (\xi^{(j)})^2 + \sum\limits_{j \geq ml +1} \lambda_j .
\end{array}
\]
\hspace*{\fill}{$\Box$}

\bigskip
Set
\begin{equation}\label{3.3}
C(l) := \sup_{k \geq 1} k^{2/l} e_k (N(0,I_l))^2.
\end{equation}
By~(\ref{1.7}), $C(l) < \infty$. For every $l \in \N$,
\begin{equation}\label{(3.4)}
e_{n_{j}} (N(0,I_l)^2 \leq n^{-2/l}_j C(l)
\end{equation}
Then one may replace the optimization problem which consists, for fixed $n$, in
minimizing the right hand side of Lemma 1 by the following optimal allocation problem:
\begin{equation}\label{optipb}
\min\{ C(l) \sum\limits^m_{j=1} \lambda_{(j-1)l+1} n^{-2/l}_j +
\sum\limits_{j \geq ml +1} \lambda_j : m, l, n_1 , \ldots , n_m \in \N, \Pi^m_{j=1} n_j \leq n \} .
\end{equation}
Set
\begin{equation}\label{mln}
m = m(n,l) := \max \{ k \geq 1 : n^{1/k} \lambda_{(k-1)l+1}^{l/2} ( \Pi^k_{j=1} \lambda_{(j-1)l+1} )^{-l/2k}
\geq 1 \} ,
\end{equation}\label{nj}
\begin{equation}\label{eq:nj}
n_j = n_j (n,l) := [n^{1/m} \lambda^{l/2}_{(j-1)l+1} ( \Pi^m_{i=1} \lambda_{(i -1)l+1} )^{-l/2m}] ,\; j \in \{ 1 , \ldots , m\},
\end{equation}
where $[x]$ denotes the integer part of $x \in \R$ and
\begin{equation}\label{ln}
l = l_n := [(\max \{ 1 , \log n \})^\vartheta ], \;\vartheta \in (0,1) .
\end{equation}
In the following theorem it is demonstrated that this
choice is at least asymptotically optimal provided the eigenvalues are regularly varying.

\begin{Thm}\label{Theorem5}
Assume the situation of Theorem 1. Consider $\hat{X}^n$ with tuning parameters defined in~(\ref{mln})-(\ref{ln}).
Then $\hat{X}^n$ is asymptotically $n$-optimal, i.e.
\end{Thm}
\[
(\E\,\| X - \hat{X}^n \|^2 )^{1/2} \sim e_n (X) \; \mbox{ as } \; n \rightarrow \infty .
\]

Note that no block quantizer with fixed block length is asymptotically optimal (see \cite{LUPA2}). As mentioned above,
$\hat{X}^n$ is not a Voronoi quantization of $X$. If $\alpha_n := \hat{X}^n(\Omega)$, then the Voronoi
quantization $\hat{X}^{\alpha_n}$ is clearly also asymptotically $n$-optimal.

\bigskip
The key property for the proof is the following $l$-asymptotics of the constants $C(l)$ defined in~(\ref{3.3}). It is
interesting to consider also the smaller constants
\begin{equation}
Q(l) := \lim_{k \to \infty} k^{2/l} e_k (\mathcal{N}(0, I_l))^2
\end{equation}
(see (\ref{1.7})).

\begin{Pro}
The sequences $(C(l))_{l \geq 1}$ and $(Q(l))_{l \geq 1}$ satisfy
\[
\lim_{l \to \infty} \frac{C(l)}{l} = \lim_{l \to \infty} \frac{Q(l)}{l} = \inf_{l \geq 1} \frac{C(l)}{l}
= \inf_{l \geq 1} \frac{Q(l)}{l} = 1.
\]
\end{Pro}
{\bf Proof.} From \cite{LUPA2} it is known that
\begin{equation}
\liminf_{l \to \infty} \frac{C(l)}{l} = 1.
\end{equation}
Furthermore, it follows immediately from~(\ref{1.7}) and~(\ref{(1.8)}) that
\begin{equation}
\lim_{l \to \infty} \frac{Q(l)}{l} = 1.
\end{equation}
(The proof of the existence of $\displaystyle \lim_{l \to \infty} C(l)/l$ we owe to S. Dereich.) For $l_0, l \in \N$ with
$l \geq l_0$, write
\[
l = n\, l_0 + m \; \mbox{with} \; n \in \N, m \in \{ 0, \ldots , l_0 -1 \} .
\]
Since for every $k \in \N$,
\[
[k^{l_0/l}]^n \; [k^{1/l}]^m \leq k ,
\]
one obtains by a block-quantizer design consisting of $n$ blocks of length $l_0$ and $m$ blocks of length 1 for
quantizing $N(0, I_l)$,
\begin{equation}\label{(3.13)}
e_k (\mathcal{N}(0,I_l))^2 \leq n e_{[k^{l_0/l}]} (\mathcal{N}(0,I_{l_0}))^2 + m e_{[k^{1/l}]} (\mathcal{N}(0,1))^2 .
\end{equation}
This implies

\begin{eqnarray*}
C(l) &  \leq & n C(l_0)  \sup_{k \geq 1} \frac{k^{2/l}}{[k^{l_0/l}]^{2/l_0}} + m C(1) \sup_{k \geq 1} \frac{k^{2/l}}{[k^{1/l}]^2}  \\
     & \leq    & 4^{1/l_0} n C(l_0) + 4 m C(1) .
\end{eqnarray*}

Consequently, using $n/l \leq 1/l_0$,
\[
\frac{C(l)}{l} \leq \frac{ 4^{1/l_0} C(l_0) }{ l_0 } + \frac{ 4 m C(1) }{ l }
\]
and hence
\[
\limsup_{l \to \infty} \frac{C(l)}{l} \leq \frac{ 4^{1/l_0} C(l_0) }{ l_0 } .
\]
This yields
\begin{equation}
\limsup_{l \to \infty} \frac{C(l)}{l} \leq \liminf_{l_0 \rightarrow \infty} \frac{C(l_0)}{l_0} = 1 .
\end{equation}
It follows from~(\ref{(3.13)}) that
\[
Q(l) \leq n Q(l_0) + m Q(1) .
\]
Consequently
\[
\frac{Q(l)}{l} \leq \frac{Q(l_0)}{l_0} + \frac{m Q(1)}{l}
\]
and therefore
\[
1 = \lim_{l \to \infty} \frac{Q(l)}{l} \leq \frac{ Q(l_0) }{ l_0 } .
\]
This implies
\begin{equation}
\inf_{l_0 \geq 1} \frac{Q(l_0)}{l_0} = 1.
\end{equation}
Since $Q(l) \leq C(l)$, the proof is complete. \hfill{$\Box$}

\bigskip
The $n$-asymptotics of  the number $m(n, l_n)l_n$ of quantized coefficients in the Karhunen-Lo\`eve expansion  in
the quantization $\hat{X}^n$ is as follows.
\begin{Lem}\label{Lemma2}
(\cite{LUPA3}, Lemma 4.8) Assume the situation of Theorem 1. Let $m(n, l_n)$ be defined by~(\ref{mln}) and~(\ref{ln}). Then
\[
m(n, l_n) l_n \sim \frac{2 \log n}{b} \; \mbox{ as } \; n \rightarrow \infty .
\]
\end{Lem}

\noindent {\bf Proof of Theorem \ref{Theorem5}.} For every $n \in \N$,
\begin{eqnarray*}
\sum\limits^m_{j=1} \lambda_{(j-1)l+1} n^{-2/l}_j
& \leq & \sum\limits^m_{j=1} \lambda_{(j-1)l+1} (n_j +1)^{-2/l} ( \frac{n_j +1}{n_j} )^{2/l} \\
& \leq & 4^{1/l} m n^{-2/ml} ( \Pi^m_{j=1} \lambda_{(j-1)l+1} )^{1/m} \\
& \leq & 4^{1/l} m \lambda_{(m-1) l+1} .
\end{eqnarray*}
Therefore, by Lemma 1 and~(\ref{(3.4)}),
\[
\E\,\| X - \hat{X}^n \|^2 \leq 4^{1/l} \frac{C(l)}{l} m l \lambda_{(m-1)l+1} + \sum\limits_{j \geq m l + 1} \lambda_j
\]
for every $n \in \N$. By Lemma~2, we have
\[
m l = m(n,l_n) l_n \sim \frac{2 \log n}{b} \; \mbox{ as } \; n \rightarrow \infty .
\]
Consequently, using regular variation at infinity with index $-b < -1$ of the function $\varphi$,
\[
m l \lambda_{(m-1)l+1} \sim m l \lambda_{ml} \sim \left( \frac{2}{b} \right)^{1-b} \psi (\log n)^{-1}
\]
and
\[
\sum\limits_{j \geq m l +1} \lambda_j \sim \frac{m l \varphi(ml)}{b-1} \sim
\frac{1}{b-1} \left( \frac{2}{b} \right)^{1-b} \psi ( \log n)^{-1} \; \mbox{ as } \; n \rightarrow \infty,
\]
where, like in Theorem 1, $\psi(x) = 1/x \varphi(x)$. Since by Proposition 1,
\[
\lim_{n \to \infty} \frac{4^{1/l_n}C(l_n)}{l_n} = 1 ,
\]
one concludes
\[
\E\,\| X - \hat{X}^n \|^2 \stackrel{<}{\sim} \left( \frac{2}{b} \right)^{1-b} \frac{b}{b-1}\, \psi ( \log n)^{-1} \; \mbox{ as }
\; n \rightarrow \infty .
\]
The assertion follows from Theorem 1. \hfill{$\Box$}

\bigskip
Let us briefly comment on the true dimension of the problem.

For $n \in \N$, let ${\cal C}_n(X)$ be the (nonempty) set of all $L^2$-optimal
$n$-quantizers. We introduce the integral number
\begin{equation}
d_n^\ast(X) := \min \left\{ \mbox{dim} \; \mbox{span} \; (\alpha)  : \alpha \in {\cal C}_n (X) \right\} .
\end{equation}
It represents the dimension at level $n$ of the functional quantization problem for $X$. Here span$(\alpha)$ denotes the linear
subspace of $H$ spanned by $\alpha$. In view of Lemma \ref{Lemma2}, a reasonable conjecture for Gaussian random vectors is $d_n^\ast(X) \sim 2 \log n/b$ in regular
cases, where $-b$ is the regularity index. We have at least the following lower
estimate in the Gaussian case.

\begin{Pro}\label{PropD}
Assume the situation of Theorem 1. Then
\[
d_n^\ast (X) \stackrel{>}{\sim} \frac{1}{b^{1/(b-1)}} \; \frac{2 \log n}{b} \; \mbox{ as } \;
n \rightarrow \infty . \vspace{0.5cm}
\]
\end{Pro}
{\bf Proof.} For every $n \in \N$, we have
\begin{equation}\label{4_2}
d_n^\ast (X) = \min \left\{ k \geq 0 : e_n ( \bigotimes^k_{j=1} N(0, \lambda_j))^2 + \sum\limits_{j \geq k+1} \lambda_j \leq e_n (X)^2 \right\}
\end{equation}
(see \cite{LUPA1}). Define
\[
c_n := \min \left\{ k \geq 0 : \sum\limits_{j \geq k + 1} \lambda_j \leq e_n (X)^2 \right\} .
\]
Clearly, $c_n$ increases to infinity as $n \rightarrow \infty$ and by (\ref{4_2}), $c_n \leq d_n^\ast (X)$ for every $n \in \N$.
Using Theorem 1 and the fact that $\psi$ is regularly varying at infinity with index $b-1$, we obtain
\[
((b-1) \psi(c_n))^{-1}  \sim  \sum\limits_{j \geq c_n +1} \lambda_j \sim e_n (X)^2
 \sim  \left( \frac{2}{b} \right)^{1-b} \frac{b}{b-1} \,\psi (\log n)^{-1}
\]
and thus
\[
\psi(c_n)  \sim \left( \frac{2}{b} \right)^{1-b} \frac{1}{b} \psi (\log n)
\sim  \psi  \left( \frac{1}{b^{1/(b-1)}} \; \frac{2 \log n}{b}\right) \; \mbox{ as } \; n \rightarrow \infty .
\]
Consequently,
\[
c_n \sim \frac{1}{b^{1/(b-1)}} \; \frac{2 \log n}{b} \; \mbox{ as } \; n \rightarrow \infty .
\]
This yields the assertion. \hfill{$\Box$}

\section{Quantizer designs and applications}\label{sec3}
\setcounter{equation}{0}
\setcounter{Assumption}{0}
\setcounter{Theorem}{0}
\setcounter{Proposition}{0}
\setcounter{Corollary}{0}
\setcounter{Lemma}{0}
\setcounter{Definition}{0}
\setcounter{Remark}{0}

In this section we are no longer interested in only
asymptotically optimal quantizers of a Gaussian process $X$,
but rather in really optimal or at least locally optimal quantizers for finite
$n\in\N$. 

As soon as the Karhunen-Lo\`eve basis $(u_j)_{j\geq 1}$ and the corresponding
eigenvalues $(\lambda_j)_{j\geq 1}$ of the Gaussian process $X$ are known,
it is possible to transform the quantization problem of $X$ in $H$ into the
quantization of $\ProdNormal$ on $l^2$ by the isometry $S: H \rightarrow l^2$
\[
	x \mapsto \left(\scalar{u_j}{x}\right)_{j\geq 1}.
\]
and its inverse
\begin{equation}
  \label{eq:isoTInv}
	S^{-1}: (l^2, \scalarOn{\cdot\;}{\cdot}{K}) \rightarrow (H,\scalar{\cdot}{\cdot}),\quad l \mapsto \sum_{j\geq 1} l_j u_j.
\end{equation}
The transformed problem then allows as we will see later on a direct access by
vector quantization methods.

The following result is straightforward.
\begin{Pro}
\label{prop:iso}
Denote by $\alpha \subset H$ an arbitrary quantizer for $X$ with associated Voronoi quantization $\widehat X^\alpha$
and Voronoi partition $\{C_a(\alpha): a\in\alpha\}$.
If
\[
        S:(H, \scalar{\cdot}{\cdot}) \rightarrow (K, \scalarOn{\cdot}{\cdot}{K})
\]
is a bijective isometry from $H \rightarrow K$, where $K$ is another separable
Hilbert space,
i.e. $S$ is linear and $\norm{Sx-Sy}_{l^2} = \norm{x-y} $ for every $x,y\in H$,
then
\begin{enumerate}
  \item $S\bigl(C_a(\alpha)\bigr) = C_{Sa}(S\alpha), \quad \text{for every } a\in \alpha$
  \item $S(\widehat X^\alpha) = \widehat{S(X)}^{S\alpha}$ is a Voronoi quantization of $S(X)$ induced by $S\alpha$
  \item $\E \min_{a\in\alpha} \norm{X-a}^2  = \E \min_{Sa\in S\alpha} \norm{S(X)-Sa}^2_{K}$
\end{enumerate}
\end{Pro}

Consequently we may focus on the quantization problem of the Gaussian random vector
\[
	\zeta := S(X)
\]
on $l^2$ with distribution
\[
	\zeta = (\zeta_j)_{j\geq 1} \sim \ProdNormal
\]
for the eigenvalues $(\lambda_j)_{j\geq 1}$ of $C_X$.
Note, that in this case $(\lambda_j)_{j\geq 1}$ also become the
eigenvalues of the covariance operator $C_\zeta$.

\subsection{Optimal Quantization of $\ProdNormal$}

Since an infinite dimensional 
quantization problem is without any modification not solvable by a finite computer
algorithm,
 we have to somehow  reduce the dimension of the problem.

Assume $\alpha$ to be an optimal $n$-quantizer for $\ProdNormal$,
then $U:=\myspan(\alpha)$ is a subspace of $l^2$ with dimension $d^\ast_n = \dim U \leq n-1$.
Consequently there exist $d^\ast_n$ orthonormal vectors in $l^2$ 
such that $\myspan(u_1, \dots, u_{d^\ast_n}) = U$.

Theorem 3.1 in \cite{LUPA1} now states,
that this orthonormal basis of $U$ can be constructed by eigenvectors of $C_\zeta$,
which correspond to the $d^\ast_n$ largest eigenvalues.
To be more precise, we get
\begin{equation}\label{eq:optBlock}
	e_n^2\biggl(\ProdNormal\biggr) = e_n^2\biggl(\bigotimes_{n=1}^{d^\ast_n} \Normal(0,\lambda_n)\biggr) + \sum_{j\geq {d^\ast_n}+1} \lambda_j.
\end{equation}
Hence it is sufficient to quantize only the finite-dimensional product measure $\bigotimes_{j=1}^{d^\ast_n} \Normal(0, \lambda_j)$
and to fill the remaining quantizer components with zeros.

Therefore we denote by $\zeta^d$ the projection of $\zeta = (\zeta_j)_{j\geq 1}$ 
on the first $d$-components, i.e. $\zeta^d = (\zeta_1, \dots, \zeta_d)$.

This approach leads for some $d\in\N$ to our first quantizer design.
\begin{quantizerDesign}[H]
	\caption{Product Quantizer for $\ProdNormal$}
	\label{des:single}
	\begin{quantizerDesignic}
		\REQUIRE Optimal $\bigotimes_{j=1}^{d}\Normal(0,\lambda_j)$-Quantizer $\alpha^{d} \subset \R^{d}$ with $\card(\alpha^{d}) \leq n$		
		\bigskip
        \QUANTIZER \[ \alpha^{\QD} := \alpha^d \times \{0\} \times \dots \]
		\QUANTIZATION \[ 
			\widehat\zeta^{\alpha^{\QD}} = \sum_{a\in\alpha^{\QD}} a \indicator{C_a(\alpha^{\QD})}{\zeta} 
			= (\widehat{\zeta^d}^{\alpha^d}, 0, \dots)
		\]
		\DISTORTION \[ 
			\E\norm{\zeta - \widehat\zeta^{\alpha^{\QD}}}^2_{l^2} 
			= e_n^2\biggl(\bigotimes_{j=1}^d \Normal(0,\lambda_j)\biggr) + \sum_{j\geq d+1} \lambda_j 
		\]
	\end{quantizerDesignic}
\end{quantizerDesign}

The claim about the distortion of $\widehat\zeta^{\alpha^{\QD}}$ becomes immediately evident from the orthogonality of the basis 
$v_j = (\delta_{ij})_{i\geq 1}$ in $l^2$ and 
\begin{align*}
	\E\norm{\zeta - \widehat\zeta^{\alpha^{\QD}}}^2_{l^2} & 
	= \E\biggnorm{\sum_{j=1}^d \Bigl(\zeta_j -
	\bigl(\widehat{\zeta^d}^{\alpha^d}\bigr)_j \Bigr) v_n + \sum_{j\geq d+1} \zeta_j v_n }^2_{l^2} \\
	& = \E \sum_{j=1}^d \Bigl( \zeta_j - \bigl(\widehat{\zeta^d}^{\alpha^d}\bigr)_j \Bigr)^2 + \sum_{j\geq d+1} \E \zeta_j^2.
\end{align*}

Unfortunately the true value of $d^\ast_n$ is only known for $n=2$, which yields $d^\ast_2 = 1$,
but from Proposition \ref{PropD} we have the lower asymptotical bound
\[
	\frac{1}{b^{1/(b-1)}} \frac{2\log n}{b} \lesssim d^\ast_n, \qquad \text{as } n\to \infty,
\]
whereas there is a conjecture for it to be $d^\ast_n \sim 2 \log n/b$.

A numerical approach for this optimal design by means of a stochastic gradient
method will be introduced in section \ref{OptiFBW},
where also some choices for the block size $d$ with regard to the quantizer size
$n$ will be given.

In addition to this direct quantization design,
we want to present some product quantizer designs for $\ProdNormal$,
which are even tractable by deterministic integration methods 
and therefore achieve a higher numerical accuracy and stationarity.
These product designs reduce furthermore the storage demand for the precomputed
quantizers when using functional quantization as cubature formulae e.g.

To proceed this way, we replace the single quantizer block $\alpha^d$ from
Quantizer
Design \ref{des:single} by the cartesian product of say $m$ smaller blocks 
with maximal dimension $l < d$.
We will refer to the dimension of these blocks also as the {\it block length}.

Let $l_i$ denote the length of the $i$-th block and set 
\[
	k_1 := 0, \qquad k_i := \sum_{\nu=1}^{i-1} l_\nu,\quad i \in \{2, \dots, m\},
\] 
then we obtain a decomposition of $\zeta^d$ into
\begin{equation}
\label{eq:blockDecomposition}
	\zeta^d = (\zeta^{(1)}, \dots, \zeta^{(m)} ), \quad \text{with} \quad 
		\zeta^{(i)} := (\zeta_{k_i+1}, \dots, \zeta_{k_i+l_i} = \zeta_{k_{i+1}}).
\end{equation}

So we state for some $l\in\N$:
\begin{quantizerDesign}[H]
	\caption{Product Quantizer for $\ProdNormal$}
	\label{des:FourierQuant}
	\begin{quantizerDesignic}
		\REQUIRE Optimal $\bigotimes_{j=k_i + 1}^{k_{i+1}}\Normal(0,\lambda_j)$-Quantizers $\ai \subset \R^{l_i}$ 
			with $\card(\ai) \leq n_i$
			for some Integers $m\in \N,\ l_1, \dots l_m \leq l,\ n_1, \dots, n_m > 1, \ \prod_{i=1}^m n_i \leq n$ solving
		\BLOCKALLOCATION \[ 
				\Biggl\{\sum_{i=1}^m e_{n_i}^2\Biggl(\bigotimes_{k=k_i + 1}^{k_{i+1}}\Normal(0,\lambda_j)\Biggr) 
				+ \sum_{j\geq k_{m+1}+1} \lambda_j \Biggr\} \rightarrow \min.
			\]
		\QUANTIZER \[ \alpha^{\QD} := \bigtimes_{i=1}^m\ai \times \{0\} \times \dots \]
		\QUANTIZATION \[ 
			\widehat\zeta^{\alpha^{\QD}} = \sum_{a\in\alpha^{\QD}} a \indicator{C_a(\alpha^{\QD})}{\zeta} 
			= (\widehat{\zeta^{(1)}}^{\alpha^{(1)}}, \dots, \widehat{\zeta^{(m)}}^{\alpha^{(m)}} , 0, \dots)
		\]
		\DISTORTION \[ 
			\E\norm{\zeta - \widehat\zeta^{\alpha^{\QD}}}^2_{l^2} 
			= \sum_{i=1}^m e_{n_i}^2\Biggl(\bigotimes_{j=k_i + 1}^{k_{i+1}}\Normal(0,\lambda_j)\Biggr) 
				+ \sum_{j\geq k_{m+1}+1} \lambda_j
		\]
	\end{quantizerDesignic}
\end{quantizerDesign}

Note that we do not use the asymptotically block allocation rules for the $n_i$ from (\ref{eq:nj}),
but perform instead the block allocation directly on the true distortion of the quantizer block and not on an estimate for them.

Next, we weaken our quantizer design,
and obtain this way the asymptotically optimal design from Theorem \ref{Theorem5}.

In fact the quantizer used for this scheme are a little bit more universal,
since they do not depend on the position of the block, but not at all more 
simply to generate.

The idea is to quantize blocks $\xi^{(i)} \sim \Normal(0,I_{l_i})$ of standard normals 
$\xi = (\xi_j)_{j\geq 1} \sim \bigotimes_{j=1}^\infty \Normal(0,1)$ 
and to weight the quantizers by
\[
	\sqrt{\lambda^{(i)}} := \biggl( \sqrt{\lambda_{k_i+1}}, \dots, \sqrt{\lambda_{k_{i+1}}} \biggr), \quad i\in\{1, \dots, m\},
\]
that is
\[
	\sqrt{\lambda^{(i)}} \otimes \ai = \Bigl\{ ( \sqrt{\lambda_{k_i+1}} a_{k_i+1}, \dots , \sqrt{\lambda_{k_{i+1}}} a_{k_{i+1}} ) : 
				a = (a_{k_i+1}, \dots, a_{k_{i+1}}) \in \ai \Bigr\}.
\]

The design for some $l \in \N$ then reads as follows:
\begin{quantizerDesign}[H]
	\caption{Product Quantizer for $\ProdNormal$}
	\label{des:NormalQuant}
	\begin{quantizerDesignic}
		\REQUIRE Optimal $\bigotimes_{j=k_i + 1}^{k_{i+1}}\Normal(0,1)$-Quantizers $\ai \subset \R^{l_i}$ 
			with $\card(\ai) \leq n_i$
			for some Integers $m\in \N,\ l_1, \dots l_m \leq l,\ n_1, \dots, n_m > 1, \ \prod_{i=1}^m n_i \leq n$ solving
		\BLOCKALLOCATION \[ 
				\Biggl\{\sum_{i=1}^m \sum_{j=k_i+1}^{k_{i+1}} \lambda_j \E \biggl(\xi_j - \Bigl(\widehat{\xi^{(i)}}^{\ai}\Bigr)_j \biggr)^2
				+ \sum_{j\geq k_{m+1}+1} \lambda_j \Biggr\} \rightarrow \min.
			\]
		\QUANTIZER \[ \alpha^{\QD} := \bigtimes_{i=1}^m \sqrt{\lambda^{(i)}} \otimes \ai \times \{0\} \times \dots \]
		\QUANTIZATION \[ 
			\widehat\zeta^{\alpha^{\QD}} = \sum_{a=(a^{(1)}, \dots, a^{(m)}, 0, \dots)\in\alpha^{\QD}} a 
				\prod_{i=1}^m \indicator{ C_{a^{(i)}}\bigl(\sqrt{\lambda^{(i)}} \otimes \alpha^{(i)} \bigr) }{\zeta^{(i)}} 
		\]
		\DISTORTION \[ 
			\E\norm{\zeta - \widehat\zeta^{\alpha^{\QD}}}^2_{l^2} 
			= \sum_{i=1}^m \sum_{j=k_i+1}^{k_{i+1}} \lambda_j \E \biggl(\xi_j - \Bigl(\widehat{\xi^{(i)}}^{\ai}\Bigr)_j \biggr)^2
				+ \sum_{j\geq k_{m+1}+1} \lambda_j
		\]
	\end{quantizerDesignic}
\end{quantizerDesign}

In the end we state explicitly the case $l=1$, 
for which the Designs \ref{des:FourierQuant} and \ref{des:NormalQuant} coincide,
and which relies only on one dimensional quantizers of the standard normal
distribution. 
These quantizers can be very easily constructed by a standard Newton-algorithm,
since the Voronoi-cells in dimension one are just simple intervals.

This special case corresponds to a direct quantization of the Karhunen-Lo\`eve expansion \eqref{eq:KL}.

We will refer to this design also as {\it scalar product quantizer}.
\begin{quantizerDesign}[H]
	\caption{Product Quantizer for $\ProdNormal$}
	\label{des:scalarQuant}
	\begin{quantizerDesignic}
		\REQUIRE Optimal $\Normal(0,1)$-Quantizers $\alpha_i \subset \R$ 
			with $\card(\alpha_i) \leq n_i$
			for some Integers $m\in \N,\ n_1, \dots, n_m > 1, \ \prod_{i=1}^m n_i \leq n$ solving
		\BLOCKALLOCATION \[ 
				\Biggl\{\sum_{j=1}^m \lambda_j e^2_{n_j}\bigl(\Normal(0,1)\bigr)
				+ \sum_{j\geq m+1} \lambda_n \Biggr\} \rightarrow \min.
			\]
		\QUANTIZER \[ \alpha^{\QD} := \bigtimes_{j=1}^m \sqrt{\lambda_j} \alpha_j \times \{0\} \times \dots \]
		\QUANTIZATION \[ 
			\widehat\zeta^{\alpha^{\QD}} = \sum_{a\in\alpha^{\QD}} a \prod_{j=1}^m \indicator{C_a(\sqrt{\lambda_j}\alpha_j)}{\zeta_j} 
			\overset{d}{=} ( \sqrt{\lambda_1} \widehat{\xi_1}^{\alpha_1}, \dots, \sqrt{\lambda_m} \widehat{\xi_m}^{\alpha_m} , 0, \dots)
		\]
		\DISTORTION \[ 
			\E\norm{\zeta - \widehat\zeta^{\alpha^{\QD}}}^2_{l^2} 
			= \sum_{j=1}^m \lambda_j e^2_{n_j}\bigl(\Normal(0,1)\bigr)
				+ \sum_{j\geq m+1} \lambda_j
		\]
	\end{quantizerDesignic}
\end{quantizerDesign}

Clearly, it follows from the decomposition (\ref{eq:blockDecomposition}),
that Design \ref{des:single} is optimal as soon the quantization of $\bigotimes_{n=1}^{d^\ast_n} \Normal(0,\lambda_n)$
is optimal.
Furthermore we obtain the proof of the asymptotically optimality for the quantizer
Designs
\ref{des:FourierQuant} and \ref{des:NormalQuant} from Theorem \ref{Theorem5}
using the tuning parameter
\begin{equation}
  \label{eq:tuningL}
	l:=l_n := [(max\{1, \log n \})^\theta ]\quad \text{ for some } \theta \in (0,1),
\end{equation}
%
i.e.
\[
	\E\norm{\zeta - \widehat\zeta^{\alpha^{\ref{des:single}}}}^2  \sim
	\E\norm{\zeta - \widehat\zeta^{\alpha^{\ref{des:FourierQuant}}}}^2_{l^2} 
	\sim \E\norm{\zeta - \widehat\zeta^{\alpha^{\ref{des:NormalQuant}}}}^2_{l^2}
	\sim \Bigl(\frac{b}{2}\Bigr)^{b-1} \frac{b}{b-1} \psi(\log n)^{-1} 
\]
as $n\rightarrow \infty$.

Using the same estimates as in the proof of Theorem \ref{Theorem5} for the
Design \ref{des:scalarQuant},
we only get
\begin{equation}\label{eq:scalarConst}
 \E\norm{\zeta - \widehat\zeta^{\alpha^{\ref{des:scalarQuant}}}}^2_{l^2}
	\lesssim \Bigl(\frac{b}{2}\Bigr)^{b-1} \frac{4C(1)(b-1)+1}{b-1} \psi(\log n)^{-1},
\end{equation}
so that we only can state, 
that Design \ref{des:scalarQuant} is rate optimal.

\begin{Remark}
Note that if we replace the assumption of optimality for the quantizer blocks by stationarity in 
Designs \ref{des:single}-\ref{des:scalarQuant},
the resulting quantizers are again stationary (but not necessary asymptotically optimal).
\end{Remark}

\subsection{Numerical optimization of quadratic functional 
quantization}\label{OptiFBW}

Optimization of the (quadratic) quantization of $\R^d$-valued random vector
has been extensively investigated since the early 1950's,
first in $1$-dimension, then in higher dimension when the cost of numerical
Monte Carlo simulation was drastically cut down
(see~\cite{GERSH}). Recent application of optimal vector quantization 
to numerics
turned out to be   much more  demanding  in terms of
accuracy. In that direction, one may cite~\cite{pagesOQ}, \cite{MRBenH}
(mainly focused on numerical optimization of the quadratic
quantization of normal distributions). To apply the methods developed in
these papers, it is more convenient to rewrite our optimization
problem with respect to the standard $d$-dimensional distribution 
${\cal N}(0,I_d)$
by simply considering the Euclidean norm derived
from the covariance matrix ${\rm Diag}(\lambda_1,\ldots, 
\lambda_{d^\ast_n})$ $i.e.$

\[
(\text{Quantizer Design \ref{des:single}})\Leftrightarrow\left\{\begin{array}{l}n\mbox{-optimal 
quantization of } \displaystyle
\bigotimes_{k=1}^{d^\ast_n} {\cal
N}(0,1)\\
\mbox{for the covariance norm }
|(z_1,\ldots,z_{d^\ast_n})|^2=
\sum_{k=1}^{d^\ast_n} \lambda_k z^2_k.
\end{array}\right.
\]

\begin{figure}
\centering
\includegraphics[width=12cm, height=8cm]{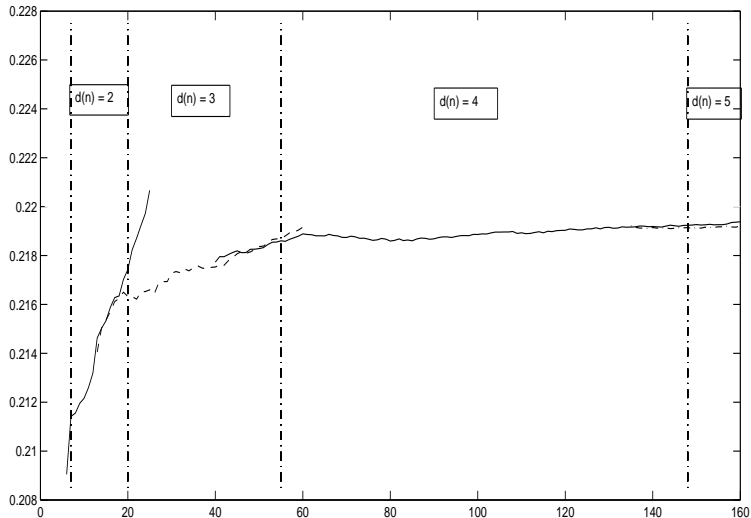}
\caption{\em  Optimal functional quantization of the Brownian motion.
$\displaystyle n\mapsto \log n \,(e_{_n}(W,L^2_{_T}))^2$,
$n\!\in\{6,\ldots,160\}$ for blocksizes $d_n \in\{2,3,4,5\}$. Vertical dashed lines: critical
dimensions for $d^\ast_n$, $e^2 \approx 7$, $e^3\approx
20$, $e^4\approx 55$, $e^5\approx 148$.}
\label{fig:dNegallogN}
\end{figure}

The main point is of course that the dimension $d^\ast_n$ is unknown.
However (see~Figure~\ref{fig:dNegallogN}), one clearly verifies on   small values
of $n$ that in the case of the Brownian Motion, i.e $b=2$ the conjecture ($d^\ast_n\sim
\log n$) is most likely true. Then  for higher
values of $n$ one relies on it to shift from one dimension to
another following the rule $d^\ast_n =d$,
$n\!\in\{e^{d},\ldots,e^{d+1}-1\}$.

\subsubsection{A toolbox for quantization optimization: a short overview}
Here is a short overview of stochastic optimization methods to compute
optimal or at least locally optimal quantizers in finite dimension. For
more details we refer to~\cite{pagesOQ} and the references therein. Let 
$Z\stackrel{d}{=} {\cal N}(0;I_d)$
and denote by $D^Z_n(x)$ the distortion function,
which is in fact the squared quantization error of a quantizer $x\!\in H^n =
(\Rd)^n$ in $n$-tuple notation, i.e.
\[
D^Z_n:  H^n \rightarrow \R, \quad x \mapsto \E \min_{1\leq i \leq n} \norm{Z-x_i}_H^2.
\]
\medskip
\noindent{\em Competitive Learning Vector Quantization ($CLVQ$).}
This procedure is a  recursive stochastic approximation gradient descent
based on the integral representation of the gradient $\nabla
D^Z_n(x),\, x\!\in H^n$ 
of the distortion as the expectation of a {\em local gradient} and a sequence of
i.i.d. random variates, i.e.
$$
\forall\, x\in H^n,\quad\nabla D^Z_n(x) = \E (\nabla
D^Z_n(x,Z))\; \text{ and }\;
Z_k\;\; i.i.d.,
\; Z_1\stackrel{d}{=} {\cal N}(0,I_d)
$$
for $\nabla D^Z_n(x) =\left(2 \int_{C_i(x)} (x_i - \xi) \Prob_Z(d\xi)
\right)_{1\leq i \leq n}$ and $\nabla D^Z_n(x, Z) =\left(2  (x_i - Z)1_{C_i(x)}(Z)
\right)_{1\leq i \leq n}$ so that, starting from  $x(0)\!\in (\R^d)^n$,
one sets
\begin{eqnarray*}
\forall\,k\ge 0,\quad x(k+1)&=&
x(k)-\frac{c}{k+1}\nabla
D^Z_n(x(k) ,Z_{k+1})
\end{eqnarray*}
where $c\!\in(0,1]$ is a real constant to be tuned. As set, this 
looks quite formal but the
operating
$CLVQ$ procedure consists of two phases at each iteration:

\smallskip
$(i)$ {\em Competitive Phase:} Search of the nearest neighbor 
$x(k)_{i*(k+1)}$ of $Z_{k+1}$
among the components of
$x(k)_i$, $i=1,\ldots,n$  (using  a ``winning convention"  in 
case of conflict on the
boundary of the Voronoi cells).

\smallskip
$(ii)$ {\em Cooperative Phase:} One moves the winning component 
toward $\zeta_{k+1}$ using a
dilatation $i.e.$
$x(k+1)_{i^*(k+1)}= {\rm Dilatation}_{\zeta_{k+1}, 
1-\frac{c}{k+1}}( x(k)_{i^*(k+1)} )$.

\smallskip This procedure is useful for small or medium values of 
$n$. For an extensive study of this procedure, which turns out to be
singular in the world of recursive stochastic approximation algorithms, we
refer to~\cite{pagesIntegration}. For general background on stochastic approximation,
we refer to~\cite{KUYI,BMP}.

\medskip
\noindent{\em The randomized ``Lloyd~I procedure".} This
is the randomization of the  stationarity based fixed point procedure 
since any optimal quantizer
satisfies 
the stationarity property:
\[
\widehat Z^{x(k+1)} = \E(Z\,|\,\widehat
Z^{x(k)}),\qquad
x(0)\subset \R^d.
\]
At every iteration the conditional expectation $\E(Z\,|\,\widehat
Z^{x(k)})$ is computed using a Monte Carlo simulation. For
more details about practical aspects of Lloyd~I procedure we refer
to~\cite{pagesOQ}. In~\cite{MRBenH}, an approach  based on  genetic 
evolutionary algorithms is developed.

For both procedures, one may substitute a sequence of quasi-random 
numbers to the usual pseudo-random sequence.
This often speeds up the rate of convergence of the method,  although 
this can only be proved (see~\cite{LAPASA}) for a
very specific class of stochastic algorithm (to which $CLVQ$ does not belong).

\smallskip  The most important step to preserve the accuracy of the 
quantization as $n$ (and $d^\ast_n$) increase is to
use the so-called {\em splitting method} which finds its origin in the proof
of the existence of an optimal $n$-quantizer: once the optimization
of a quantization grid of size $n$ is achieved, one specifies the
starting grid for the size $n+1$ or more generally $n+ \nu$, $\nu\ge 1$, by
merging the optimized grid of size $n$ resulting from the former
procedure with  $\nu$ points sampled independently from
the normal distribution with probability density proportional to
$\varphi^{\frac{d}{d+2}}$ where $\varphi$ denotes the p.d.f. of
${\cal N}(0;I_d)$. This rather  unexpected choice is motivated by
the  fact that this distribution provides  the
lowest    {\em in average}  random quantization error  (see~\cite{COH}).

\smallskip
  As a result, to be downloaded on the website~\cite{Website} devoted 
to quantization:
\medskip

\centerline{\tt www.quantize.maths-fi.com}

\medskip $\circ$   {\em Optimized  stationary codebooks for $W$}: in
practice, the $n$-quantizers $\alpha:=\alpha^{d^\ast_n}$ of the distribution
$\otimes_{k=1}^{d^\ast_n}{\cal N}(0,\lambda_k)$, $n\!=\!1$ up to
$10\,000$ ($d^\ast_n$ runs from $1$ up to $9$).

\smallskip $\circ$   {\em Companion  parameters}:

\smallskip
$\quad$ -- distribution of $\widehat W^{\gamma}$: $\P(\widehat
W^{\gamma}\!=x_i)=\P(\widehat Z_{d^\ast_n}^{\alpha}\!= \alpha_i)$.

$\quad$ -- The quadratic quantization error: $\|W-\widehat 
W^{\gamma^n}\|_{L^2_T}$.
%

\subsection{Application to the Brownian motion on $L^2([0,T], dt)$}\label{sec3.2}
We present in this subsection numerical results for the above quantizer designs 
applied to the Brownian motion $W$ on the Hilbert space $\bigl(L^2([0,T],dt), \norm{\cdot}_{L^2_T} \bigr)$.

Recall that the eigenvalues of $C_W$ read
\begin{equation*}
  \label{eq:eigenvalues}
	\lambda_j = \left(\frac{T}{\pi(j-1/2)}\right)^2, \quad j\geq 1
\end{equation*}
and the eigenvectors
\[
	u_j  = \sqrt{\frac{2}{T}} \sin(t/\sqrt{\lambda_j}), \quad j\geq 1
\]
which imply a regularity index of $b=2$ for the regularly varying function
\begin{equation*}
  \label{eq:regVarW}
	\varphi(x) := \left(\frac{T}{\pi}\right)^2 x^{-2}.
\end{equation*}

%
%

Let $\alpha$ be a quantizer for $\ProdNormal$, then for $S^{-1}$ from 
\eqref{eq:isoTInv}
\begin{equation}
  \label{eq:productGauss}
	\gamma := S^{-1}\alpha = \Bigl\{ t\mapsto \sqrt{\frac{2}{T}} \sum_{j\geq 1} a_j \sin\bigl(\pi(j-1/2)t/T\bigr): (a_1, a_2, \dots ) \in \alpha \Bigr\}
\end{equation}
provides a quantizer for $W$, 
which produces the same quantization error as $\alpha$ and is stationary iff
$\alpha$ is. 
Furthermore we can restrict w.l.o.g. to the case $T=1$.


Concerning the numerical construction of a quantizer for the Brownian motion we need access to precomputed stationary quantizers of
$\bigotimes_{j=k_i+1}^{k_{i+1}} \Normal(0,\lambda_j)$ and $\bigotimes_{j=k_i+1}^{k_{i+1}} \Normal(0,1)$ 
for all possible combinations of the block allocation problem.
As soon as these quantizers are computed,
we can perform the Block Allocation of the quantizer Designs to produce optimal
Quantizers for $\ProdNormal$.

For the quantizers of Design \ref{des:single} we used the stochastic algorithm
from section \ref{OptiFBW}, 
whereas for Designs \ref{des:FourierQuant} - \ref{des:scalarQuant} we could
employ deterministic procedures for the integration on the Voronoi cells with
max. block lengths $l=2$ respectively $l=3$,
which provide a maximum level of stationarity, i.e. $\norm{\nabla D_n} \leq 1e^{-8}$.

\begin{table}[H]
\centering
\begin{tabular}{*{2}{r|}c}
	\multicolumn{1}{c|}{n}  &	\multicolumn{1}{c|}{$d_n$} & $\E\norm{W - \widehat W^{\gamma^{\ref{des:single}}}}^2_{L^2_T}$ \\ \hline
	 1 &  1 & 0.5000\\
	 5 &  1 & 0.1271\\
	 10 & 2 & 0.0921 \\
	 50 & 3 & 0.0558 \\
	 100 & 4 &  0.0475\\
	 500  & 6 & 0.0353 \\
	 1000 & 6 & 0.0318\\
	 5000  & 8 & 0.0258\\
	 10000 & 9 & 0.0238\\
\end{tabular}
\caption{Quantizer Design \ref{des:single}}
\label{tab:single}
\end{table}

\begin{table}[H]
\centering
\begin{tabular}{*{3}{r|}c}
	\multicolumn{1}{c|}{n}  &	\multicolumn{1}{c|}{$n_i$} & \multicolumn{1}{c|}{$l_i$} & $\E\norm{W - \widehat W^{\gamma^{\ref{des:FourierQuant}}}}^2_{L^2_T}$ \\ \hline
	 1 & 1 & 1 & 0.5000\\
	 5 & 5 & 1 & 0.1271\\
	 10 & 10 & 1 & 0.0921 \\
	 50 & $25 \times 2 = 50$ & $ 2+1 = 3$ & 0.0580 \\
	 100 & $50 \times 2 = 50$ & $2+1 = 3$ & 0.0492 \\
	 500 & $100 \times 2 = 500$ & $2+1 = 3$ & 0.0372 \\
	 1000 & $111 \times 3 \times 3 = 999$ & $2 + 1 + 2 = 5$ & 0.0339\\
	 5000 & $166 \times 10 \times 3 = 4980$ & $2 + 2 + 2 = 6$ & 0.0276\\
	 10000 & $208 \times 12 \times 4 = 9984$ & $2 + 2 + 2 = 6$ & 0.0255\\
	 100000 & $277 \times 20 \times 6 \times 3 = 99720$ & $2 + 2 + 2 + 2 = 8$ & 0.0206\\
\end{tabular}
\caption{Quantizer Design \ref{des:FourierQuant}, $l=2$}
\label{tab:Fourier}
\end{table}

\begin{table}[H]
\centering
\begin{tabular}{*{3}{r|}c}
	\multicolumn{1}{c|}{n}  &	\multicolumn{1}{c|}{$n_i$} & \multicolumn{1}{c|}{$l_i$} & $\E\norm{W - \widehat W^{\gamma^{\ref{des:NormalQuant}}}}^2_{L^2_T}$ \\ \hline
	 1 & 1 & 1 & 0.5000\\
	 5 & 5 & 1 & 0.1271\\
	 10 & $5 \times 2$ & $1+1=2$ & 0.0984 \\
	 50 & $10 \times 5 = 50$ & $ 1+2 = 3$ & 0.0616 \\
	 100 & $12 \times 4 \times 2 = 96$ & $1+1+1 = 3$ & 0.0513 \\
	 500 & $16 \times 5 \times 3 \times 2 = 480$ & $1+1+1+1 = 4$ & 0.0387 \\
	 1000 & $20 \times 25 \times 2 = 1000$ & $1 + 2 + 1 = 5$ & 0.0350\\
	 5000 & $26 \times 8 \times 8 \times 3 = 4992$ & $1 + 1 + 2 + 2 = 6$ & 0.0285\\
	 10000 & $25 \times 36 \times 11 = 9900$ & $1 + 2 + 3 = 6$ & 0.0264\\
	 100000 & $33 \times 55 \times 11 \times 5 = 99825$ & $1 + 2 + 2 + 3 = 8$ & 0.0211\\
\end{tabular}
\caption{Quantizer Design \ref{des:NormalQuant}, $l=3$}
\label{tab:Normal}
\end{table}

\begin{table}[H]
\centering
\begin{tabular}{*{3}{r|}c}
	\multicolumn{1}{c|}{n}  &	\multicolumn{1}{c|}{$n_i$} & \multicolumn{1}{c|}{$m$} & $\E\norm{W - \widehat W^{\gamma^{\ref{des:scalarQuant}}}}^2_{L^2_T}$ \\ \hline
	 1 & 1 & 1 & 0.5000\\
	 5 & 5 & 1 & 0.1271\\
	 10 & $5 \times 2$ & $2$ & 0.0984 \\
	 50 & $12 \times 4 = 48$ & $ 2$ & 0.0616 \\
	 100 & $12 \times 4 \times 2 = 96$ & $3$ & 0.0513 \\
	 500 & $16 \times 5 \times 3 \times 2 = 480$ & $4$ & 0.0387 \\
	 1000 & $23 \times 7 \times 3 \times 2 = 966$ & $4$ & 0.0352\\
	 5000 & $26 \times 8 \times 4 \times 3 \times 2 = 4992$ & $5$ & 0.0286\\
	 10000 & $26 \times 8 \times 4 \times 3 \times 2 \times 2 = 9984$ & $6$ & 0.0264\\
	 100000 & $34 \times 10 \times 6 \times 4 \times 3 \times 2 \times 2 = 97920$ & $7$ & 0.0213\\
\end{tabular}
\caption{Quantizer Design \ref{des:scalarQuant}}
\label{tab:scalar}
\end{table}

\begin{figure}[htbp]
	\centering
	\includegraphics[width=0.9\textwidth]{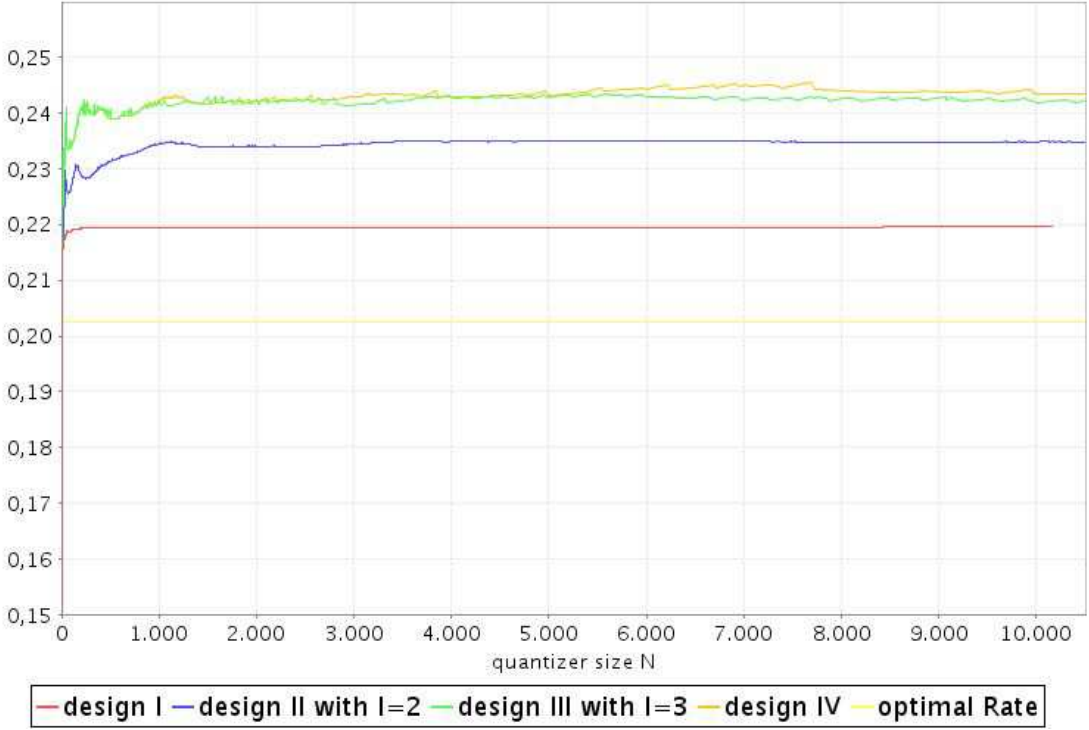}
	\caption{Asymptotics for $n\mapsto \log n \times \E\norm{W-\widehat
	W^\gamma}^2_{L^2}$ for the Designs I, II, III and IV.}
	\label{fig:quantizerPerf}
\end{figure}

The asymptotical performance of the quantizer designs in view of Theorem \ref{Theorem5}, i.e.
\[
	n \mapsto \log n\ \E\norm{W-\widehat W^\gamma}^2_{L_T^2}.
\]
is presented in Figure \ref{fig:quantizerPerf},
where the quantization coefficient is evaluated for the Brownian Motion on $[0,1]$ with 
$\varphi(j) = \pi^{-2} j^{.2}$ as
\[
	\Bigl(\frac{b}{2}\Bigr)^{b-1} \frac{b}{b-1} \pi^{-2} = \frac{2}{\pi^2} \approx 0.20264237... \ .
\]

As expected we have
\[
	\E\norm{W-\widehat W^{\gamma^{\ref{des:FourierQuant}}}}^2_{L^2_T} \leq \E\norm{W-\widehat W^{\gamma^{\ref{des:NormalQuant}}}}^2_{L^2_T} 
		\leq \E\norm{W-\widehat W^{\gamma^{\ref{des:scalarQuant}}}}^2_{L^2_T}
\]
and by (\ref{eq:scalarConst}),
\[
	\log n\ \E\norm{W-\widehat W^{\gamma^{\ref{des:scalarQuant}}}}^2_{L^2_T} \lesim
	\frac{4C(1)+1}{2} \log n\ e_n(W)^2 = 5.9414\ldots \log n\ e_n(W)^2 \sim 1.2040\ldots
\]
assuming $C(1) = Q(1)$.

Although the Designs \ref{des:single}, \ref{des:FourierQuant} and
\ref{des:NormalQuant} are asymptotically equivalent,
we can observe a great superiority of Designs \ref{des:single} and 
\ref{des:FourierQuant} compared to Design \ref{des:NormalQuant}.

\begin{figure}[htbp]
  \centering
  \begin{minipage}[c]{0.33\textwidth}
    \includegraphics[width=\textwidth]{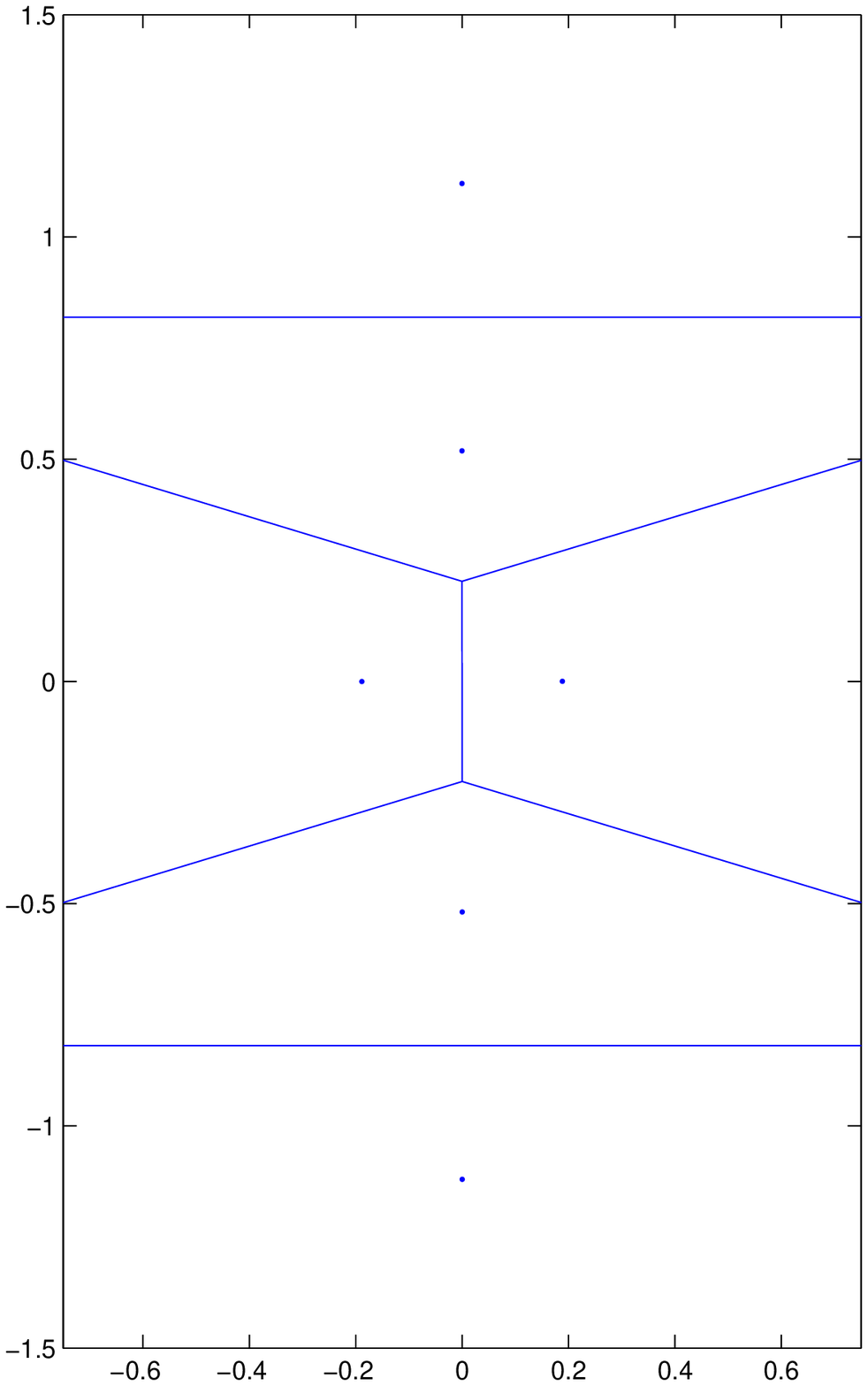}
    \end{minipage}%
  \begin{minipage}[c]{.33\textwidth}
    \includegraphics[width=\textwidth]{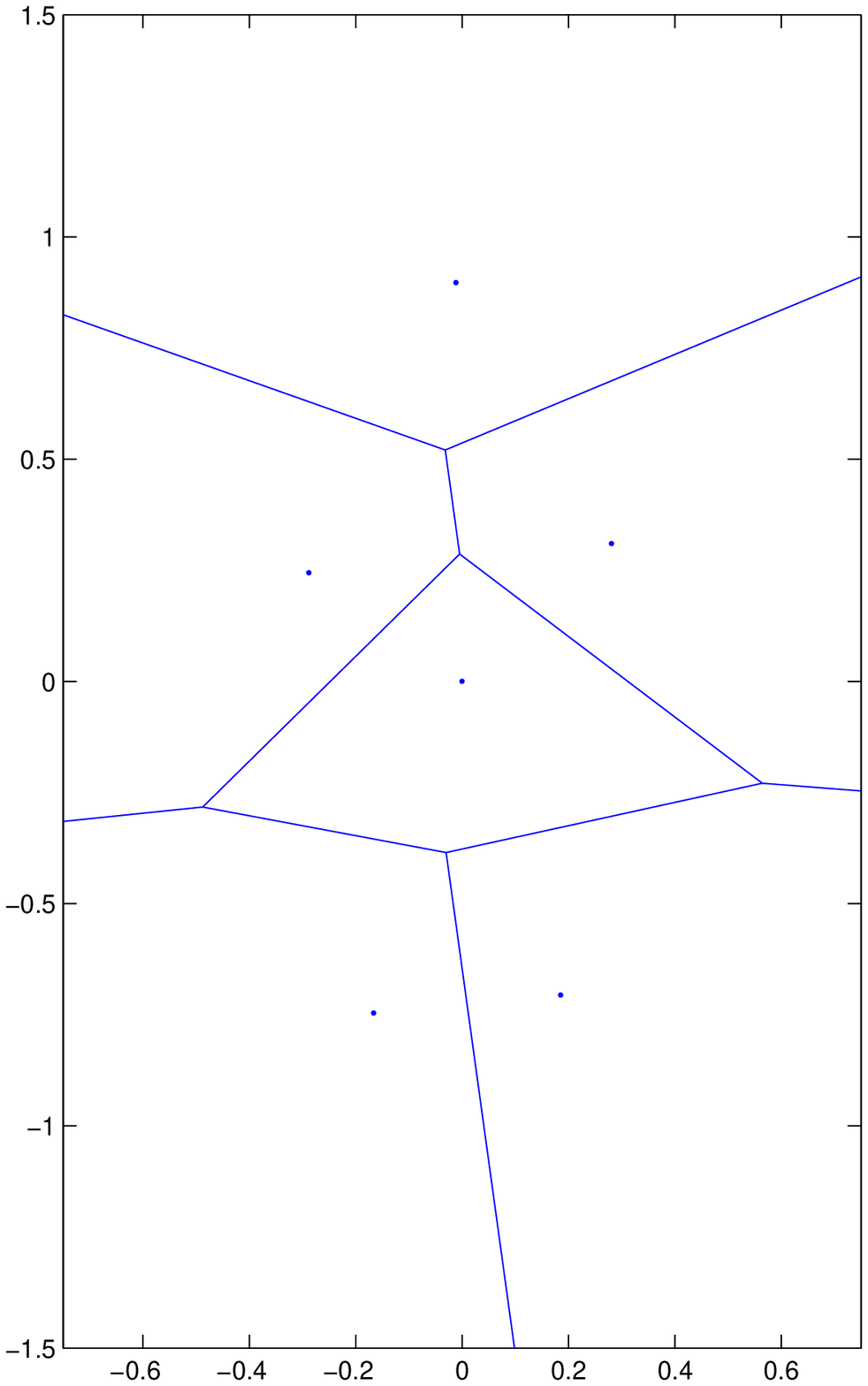}
  \end{minipage}
  \begin{minipage}[c]{0.33\textwidth}
    \includegraphics[width=\textwidth]{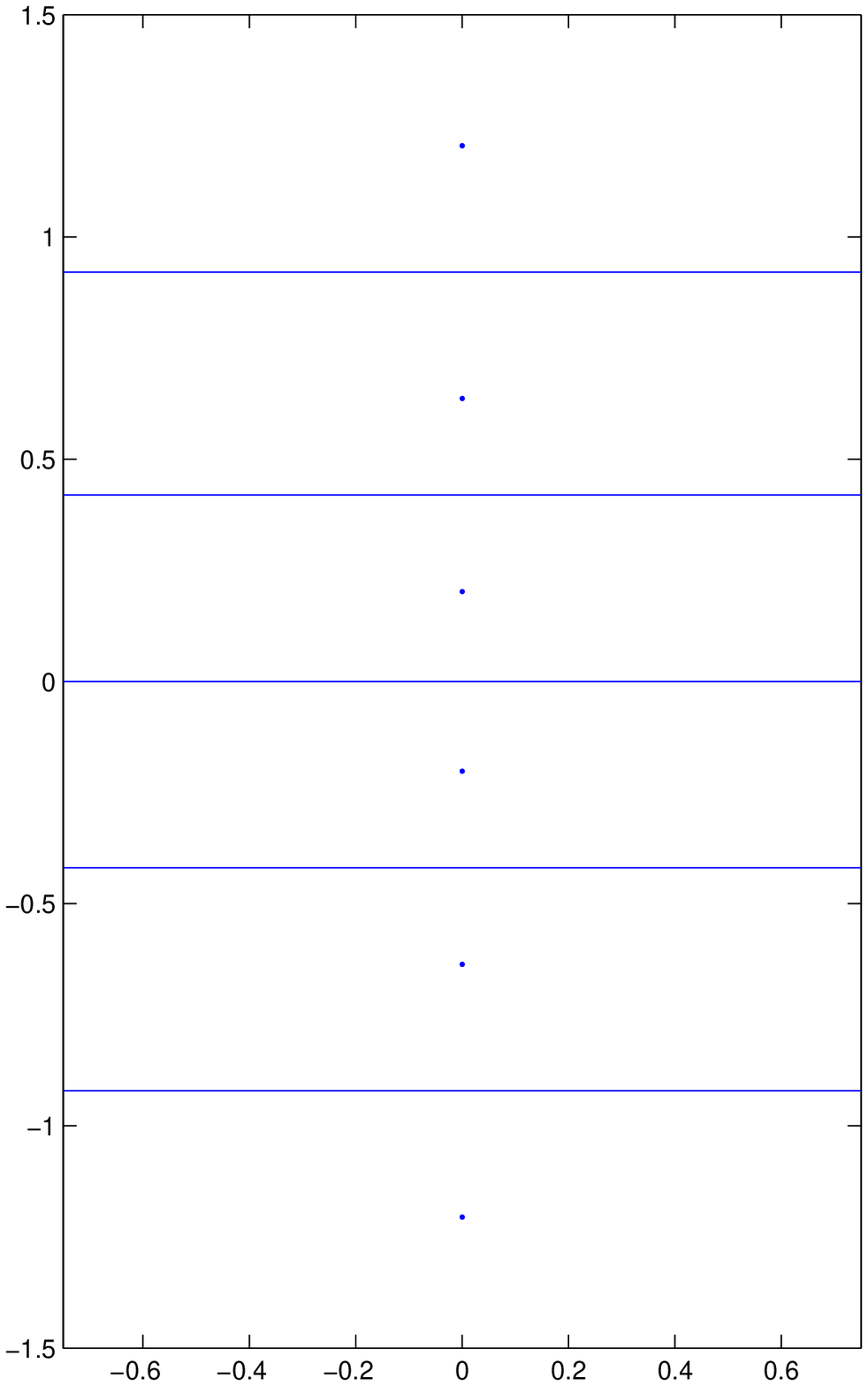}
   \end{minipage}%
	\caption{Quantizers of size $n=6$ for $\ProdNormal$ generated by Designs II, III and IV (from left to right) 
		and projected on the eigenspace corresponding to $\lambda_1$ and $\lambda_2$}
	\label{fig:quantizerComp}
\end{figure}

This is mainly caused by the better adaption to the rapidly decreasing sequence of the eigenvalues.
To give an impression of this geometrical superior adaption, we illustrate the case $n=6$ in Figure \ref{fig:quantizerComp}.
The quantizers for $\ProdNormal$ in the figure are projected onto the first two dimensions.
Within that subspace, 
quantizer \ref{des:scalarQuant} is a product quantizer of $\alpha^1\times\{0\}$,
hence the rectangular shape of the Voronoi cells.

As quantizer \ref{des:NormalQuant} was formerly optimized for the symmetrically
distribution $\Normal(0,I_2)$,
there are still to many points in the subspace generated by the eigenvector of $\lambda_2$,
which cannot be accomplished by the weightening tensor product $\sqrt{\lambda^{(i)}} \otimes \ai$.

Concerning quantizer \ref{des:FourierQuant}, we see the possibly best quantizer at level 6 for $\ProdNormal$, 
since the quantizer Design \ref{des:FourierQuant} produces the same quantizer
for $N=6$ regardless of $l=2$ or $l=3$ and is therefore equivalent to Design \ref{des:single}.


\begin{figure}[htbp]
  \centering
    \includegraphics[width=0.8\textwidth]{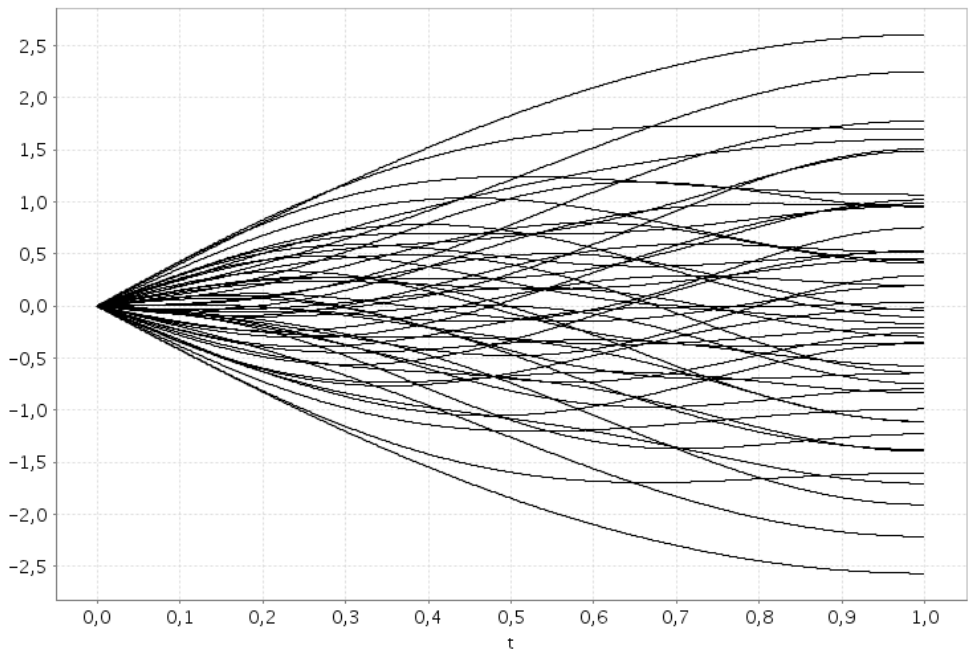}
 	\caption{A stationary quantizer for $W$ on $L^2([0,1], dt)$ generated by
 	Design \ref{des:single}, size $n=50$ and $d^\ast_n=3$ }
	\label{fig:BMquant}
\end{figure}

\begin{figure}[htbp]
  \centering
    \includegraphics[width=0.75\textwidth]{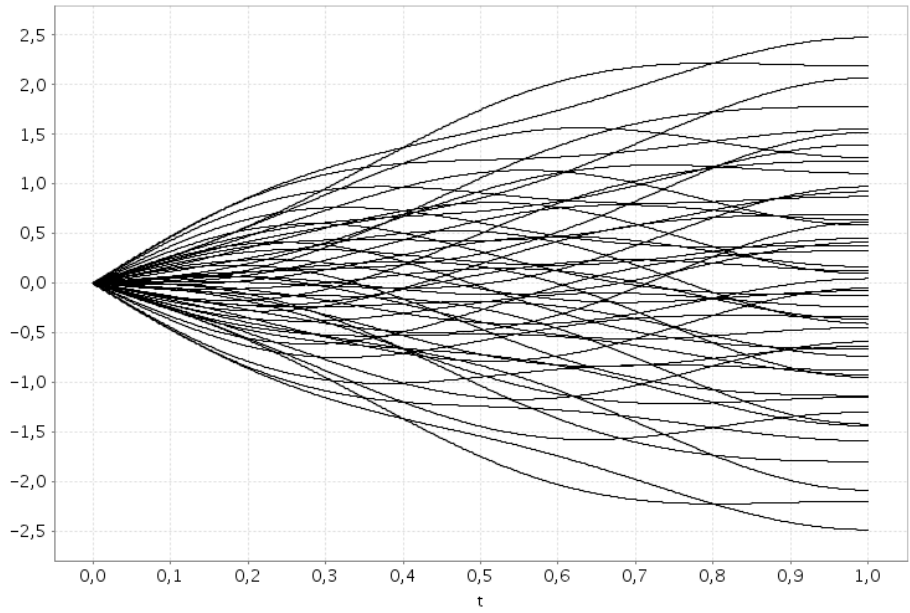}
    \includegraphics[width=0.75\textwidth]{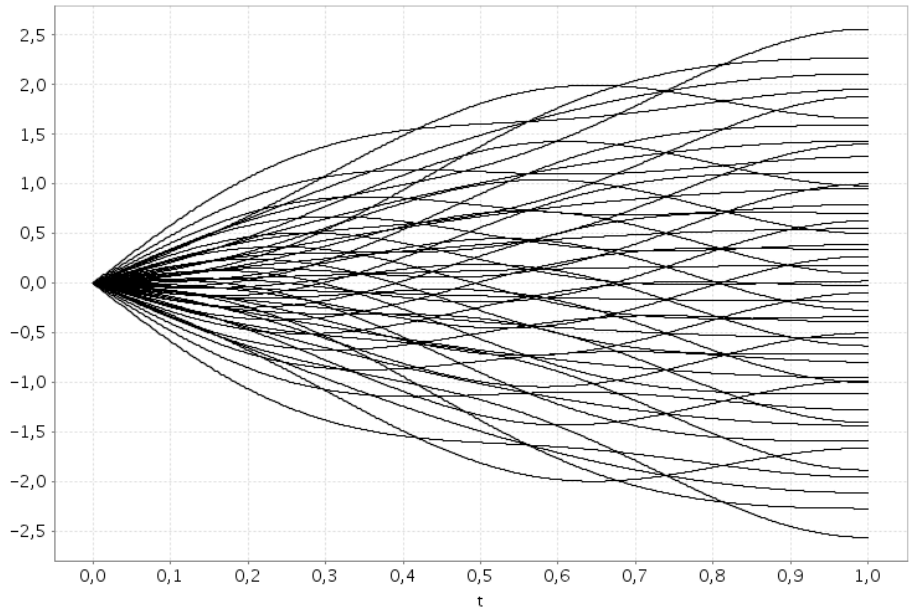}
    \includegraphics[width=0.75\textwidth]{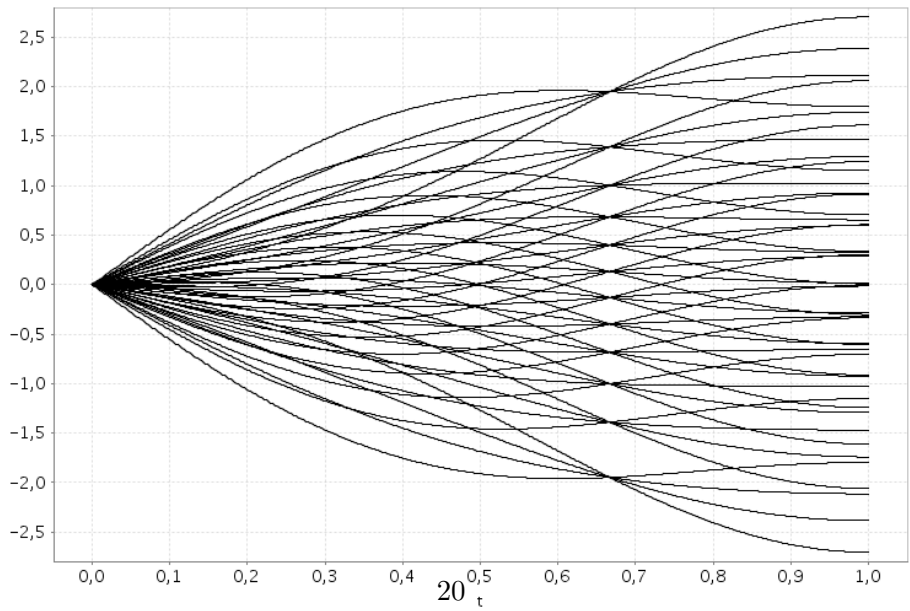}
 	\caption{Stationary quantizers of size $n=50$ for $W$ on $L^2([0,1], dt)$ generated by Designs \ref{des:FourierQuant}
 		- \ref{des:scalarQuant} (from top to bottom). See also tables \ref{tab:Fourier} - \ref{tab:scalar}.}
	\label{fig:BMquant2}
\end{figure}


\subsection{Application to Riemann-Liouville processes}

We consider Riemann-Liouville processes in
$H = L^2 ([0,T], dt)$. For $\rho \in (0, \infty)$, the
Riemann-Liouville process $X^\rho = (X^\rho_t)_{t \in [0,T]}$ on $[0,T]$ is defined by
\begin{equation}\label{(1.13)}
X^\rho_t := \int^t_0 (t - s)^{\rho- \frac{1}{2}} d W_s
\end{equation}
where $W$ is a standard Brownian motion.

Its
covariance function is given by
\begin{equation}\label{(2.1)}
\E\,X^\rho_s X^\rho_t = \int^{s \wedge t}_0 (t - r)^{\rho- \frac{1}{2}} (s-r)^{\rho - \frac{1}{2}} dr .
\end{equation}
Using $\rho \wedge \frac{1}{2}$-H{\"o}lder continuity of the application
$t \mapsto X^\rho_t$ from [0,T] into $L^2(\P)$ and
the Kolmorogov criterion one checks that $X^\rho$ has a pathwise continuous modification so that we may assume without
loss of generality that $X^\rho$ is pathwise continuous. In particular, $X^\rho$ can be seen as
a centered Gaussian random vector with values in
\[
H = L^2 ([0,T], dt).
\]
The following high-resolution formula is a consequence of a theorem by Vu and
Gorenflo \cite{tuan} on singular values of Riemann-Liouville integral operators
\begin{equation}
R_\beta \,g(t) = \frac{1}{\Gamma(\beta)} \int^t_0 (t-s)^{\beta-1} g(s) ds , \hskip 1 cm \beta \in (0, \infty).
\end{equation}

For every $\rho \in (0, \infty)$,
\begin{equation}\label{3_3}
e_n (X^\rho) \sim T^{\rho+1/2} \pi^{-(\rho + \frac{1}{2})} (\rho + 1/2)^\rho (
\frac{2 \rho + 1 }{2 \rho} )^{1/2} \Gamma( \rho + 1/2) (\log n)^{-\rho} \;
\mbox{ as } \; n \rightarrow \infty .
\end{equation}

This can be seen as follows. For $\beta > 1/2$, the Riemann-Liouville fractional
integral operator $R_\beta$ is a bounded operator from
$L^2 ([0,T], dt)$ into $ L^2 ([0,T], dt)$. The covariance operator
\[
C_\rho : L^2 ([0,T], dt) \rightarrow L^2 ([0,T], dt)
\]
of $X^\rho$ is given by the Fredholm transformation
\[
C_\rho g(t) = \int^T_0 g(s) \E X^\rho_s X^\rho_t ds .
\]
Using~(\ref{(2.1)}), one checks that $C_\rho$ admits a factorization
\[
C_\rho = S_\rho S^{*}_\rho ,
\]
where
\[
S_\rho = \Gamma (\rho +1/2 ) R_{ \rho + \frac{1}{2}} .
\]
Consequently, it follows from Theorem 1 in \cite{tuan} that the eigenvalues $\lambda_1 \geq \lambda_2 \geq \ldots > 0$ of $C_\rho$ satisfy
\begin{equation}\label{(2.3)}
\lambda_j \sim T^{2\rho+1} \Gamma ( \rho +1/2)^2 ( \pi j)^{-(2 \rho + 1)} \;
\mbox{ as } \; j \rightarrow \infty.
\end{equation}
Now (\ref{3_3}) follows from Theorem 1 (with $\varphi(x) = T^{2\rho+1} \Gamma(\rho + 1/2)^2 \pi^{-b} x^{-b}$ and $b = 2 \rho + 1)$.

\bigskip
An immediate consequence for fractionally integrated Brownian motions on $[0,T]$ defined by
\begin{equation}
Y^\beta_t := \frac{1}{\Gamma(\beta)} \int^t_0 (t-s)^{\beta-1} W_s ds
\end{equation}
for $\beta \in (0, \infty)$ is as follows.

For every $\beta \in (0, \infty),$
\[
e_n (Y^\beta) \sim T^{\beta+1} \pi^{-(\beta+1)} (\beta+1)^{\beta + \frac{1}{2}}
 ( \frac{2 \beta + 2}{2 \beta +1})^{1/2} (\log n)^{-(\beta + \frac{1}{2})} \; \mbox{ as } \; n \to \infty .
\]

In fact, for $\rho > 1/2$, the Ito formula yields
\[
X^\rho_t = (\rho - \frac{1}{2} ) \int^t_0 (t-s)^{\rho - \frac{3}{2}} W_s ds .
\]
Consequently,
\[
Y^\beta_t = \frac{1}{\beta \Gamma (\beta)} \beta \int^t_0 (t-s)^{\beta + \frac{1}{2} - \frac{3}{2} }
W_s ds = \frac{1}{\Gamma (1+\beta)} X_t^{\beta + \frac{1}{2} } .
\]
The assertion follows.

\bigskip One  further consequence is a precise relationship between the quantization errors of Riemann-Liouville processes
and fractional Brownian motions. The fractional Brownian motion with Hurst exponent $\rho \in (0,1]$
is a centered pathwise continuous Gaussian process
$Z^\rho = (Z^\rho_t)_{t \in [0,T]}$ having the covariance function
\begin{equation}
\E\,Z^\rho_s Z^\rho_t = \frac{1}{2} ( s^{2 \rho} + t^{2 \rho} - \mid s-t \mid^{2 \rho} ) .
\end{equation}

For every $\rho \in (0,1)$,
\begin{equation}\label{3_7}
e_n (X^\rho) \sim \frac{\Gamma(\rho + 1/2 )}{ ( \Gamma ( 2 \rho + 1) \sin ( \pi \rho))^{1/2} }
e_n (Z^\rho) \; \mbox{ as } \; n \rightarrow \infty .
\end{equation}
\bigskip
In fact, by \cite{LUPA2}, we have
\[
e_n (Z^\rho) \sim T^{\rho+1/2} \pi^{-(\rho + \frac{1}{2})} ( \rho + 1/2)^\rho
\left( \frac{2 \rho +1}{2 \rho}\right)^{1/2} ( \Gamma (2 \rho +1)
\sin (\pi \rho))^{1/2} ( \log n)^{-\rho} , n \rightarrow \infty .
\]
Combining this formula with (\ref{3_3}) yields the assertion (\ref{3_7})

\bigskip
Observe that strong equivalence $e_n (X^\rho) \sim e_n (Z^\rho)$ as $n \rightarrow \infty$
is true for exactly two values of $\rho \in (0,1)$, namely for $\rho = 1/2$ where even $e_n (X^{1/2} ) = e_n (Z^{1/2} ) = e_n (W)$ and, a bit
mysterious, for $\rho = 0.81557\ldots$

The basic example (among Riemann-Liouville processes) is $X^{1/2} = W$ and $H = L^2 ([0,T], dt)$, where
\begin{equation}
\lambda_j = T^2( \pi(j- 1/2))^{-2} , \;u_j(t) = \sqrt{\frac{2}{T}} \; \mbox{sin} \,
\left(t/\sqrt{\lambda_j} \right),\;  j \geq 1
\end{equation}
(see Section \ref{sec3.2}).

Since for $\delta, \rho \in (0, \infty)$,
\[
X^{\delta + \rho} = \frac{\Gamma(\delta + \rho + \frac{1}{2})}{\Gamma( \rho + \frac{1}{2})} R_\delta (X^\rho) ,
\]
one gets expansions of $X^{\delta + \rho}$ from Karhunen-Lo\`eve expansions of $X^\rho$. In particular,
\[
X^{\delta + \frac{1}{2}} = \Gamma (\delta+1) \sum\limits^\infty_{j=1} \sqrt{\lambda_j} \xi_j R_\delta (u_j).
\]
However, the functions $R_\delta(u_j), j \geq 1$, are not orthogonal in $H$ so that the nonzero
correlation between the components of $(\xi^{(j)} - \widehat{\xi^{(j)}})$ prevents the previous estimates for $\E \|X-\widehat X^n\|^2$  given in
Lemma~1 from working in this setting in the general case.

However, when $l=1$ (scalar product quantizers made up with blocks of fixed length
$l=1$, see Design \ref{des:single}), one checks that these estimates still stand as equalities
since
orthogonality can now be substituted by the independence of $\xi_j - \hat{\xi}_j$
and stationarity property~(\ref{(3.1)}) of the quantizations $\hat{\xi}_j, j \geq 1$.
It is often good enough for applications
to use  scalar product quantizers (see
\cite{LUPA1},
\cite{FQFinance}). If, for instance $\delta = 1$, then
\[
X := X^{3/2} = \sum\limits^\infty_{j=1} \sqrt{\lambda_j} \xi_j R_1 (u_j)  ,
\]
where
\[
R_1 (u_j) (t) =\sqrt{\frac{2 \lambda_j}{T}} (1 - \cos(t / \sqrt{\lambda_j})).
\]
Note that  $\displaystyle \| R_1 (u_j) \|^2 = T^2 \mu_j (3 - 4
(-1)^{j-1}\sqrt{\mu_j} ),\; j \geq 1$, where $\lambda_j = T^2 \mu_j$. Set
\[
\hat{X}^n = \sum\limits^m_{j=1} \sqrt{\lambda_j} \hat{\xi}_j R_1 (u_j) .
\]
The quantization $\widehat X^n$ is non Voronoi  (it is
related to the Voronoi tessellation  of $W$) and satisfies
\begin{equation}
\E\|X-\widehat X^n\|^2= \sum^m_{j = 1}  T^4\mu^2_j (3 - 4 (-1)^{j-1}\sqrt{\mu_j} )
e_{n_j} (N(0, 1))^2 + \sum_{j \geq m +1} T^4 \mu^2_j (3-4(-1)^{j-1}
\sqrt{\mu_j}).
\end{equation}
It is possible to optimize the (scalar product) quantization error using this expression instead of~(\ref{optipb}).
As concerns asymptotics, if the  parameters  are tuned following~(\ref{mln})-(\ref{ln}) with $l=1$ and
$\lambda_j$ replaced by
\[
\nu_j := T^4 \mu^2_j (3 + 4 \sqrt{\mu_j} ) \sim 3 \pi^{-4} j^{-4}
\quad\mbox{ as }\quad n\to \infty, \]
and using (\ref{3_3}) gives
\begin{equation}
( \E\,\| X - \hat{X}^n \|^2)^{1/2} \stackrel{<}{\sim} \left( \frac{3(12 C(1)+1)}{4}\right)^{1/2} e_n (X) \; \mbox{ as } \; n \rightarrow
\infty .
\end{equation}
Numerical experiments seem to confirm that $C(1) = Q(1)$. Since $Q(1) = \pi \sqrt{3}/2$ (see \cite{Foundations}, p. 124), the
above upper bound is then
\[
\left( \frac{3(6 \pi \sqrt{3} +1)}{4}\right)^{1/2} = 5.02357 \ldots
\]

\bibliography{literatur}

\end{document}